\documentclass[conference]{IEEEtran}

\ifCLASSINFOpdf
\else
\fi

\usepackage{a4wide,amsfonts,amsmath,latexsym,amssymb,euscript,graphicx,units,mathrsfs,hyperref,amsthm,
xcolor}

\newcommand{\R}{\mathbb{R}}
\newtheorem{prop}{Proposition}[section]
\newtheorem{cor}[prop]{Corollary}
\newtheorem{lemma}[prop]{Lemma}
\newtheorem{remark}[prop]{Remark}

\newtheorem{theorem}[prop]{Theorem}
\newtheorem{definition}[prop]{Definition}

\begin{document}

\title{Integration by parts and \\
representation of information functionals}

\author{\IEEEauthorblockN{Ivan Nourdin}
\IEEEauthorblockA{Luxembourg University\\
Luxembourg City, Luxembourg\\
Email: ivan.nourdin@uni.lu}
\and
\IEEEauthorblockN{Giovanni Peccati}
\IEEEauthorblockA{Luxembourg University\\
Luxembourg City, Luxembourg\\
Email: giovanni.peccati@uni.lu}
\and
\IEEEauthorblockN{Yvik Swan}
\IEEEauthorblockA{Universit\'e de Li\`ege\\
Li\`ege, Belgium\\
Email: yswan@ulg.ac.be}}

\maketitle

\begin{abstract}
  We introduce a new formalism for computing expectations of
  functionals of arbitrary random vectors, by using generalised
  integration by parts formulae. In doing so we extend recent
  representation formulae for the score function introduced in
  \cite{nourdin2013entropy} and also provide a new proof of a central
  identity first discovered in \cite{guo2005mutual}. We derive a
  representation for the standardised Fisher information of sums of
  i.i.d. random vectors which we use to provide rates of convergence in information theoretic central
  limit theorems (both in Fisher information distance and in relative
  entropy) and a  Stein bound for Fisher
  information distance.
\end{abstract}

\IEEEpeerreviewmaketitle

\section{Introduction}
\label{sec:foreword}

Let $X$ be a random vector in $\R^d$, with differentiable density
$f$. The \emph{score function} $\rho_X(x) = \nabla \log f(x)$ 
has
long been known to provide useful 
handles on the law of $X$. 
A much less studied object is the \emph{Stein matrix} of $X$, defined
in (\ref{eq:steinmatrix}) which
can be interpreted as a counterpart to the score where,
rather than taking log-derivatives, one considers a special form of
integration.  This matrix (whose properties when $d=1$ are closely
related to the so-called \emph{$w$-function}, see \cite{CP89})
has only recently started to attract the attention of the community,
see e.g. \cite{NPReveillac,nourdin2013entropy,LNP14}. We refer the
reader to \cite{airault2010stokes} for a detailed
study.  

In this paper we explore the connexion between the score and the Stein
 matrix of an arbitrary random vector $X$.
Rather than defining these two quantities explicitly in terms of the
density $f$, we choose to characterise them by their behaviour through
specialised integration by parts formulae 
(see equations (\ref{eq:scorefunction}) and
(\ref{eq:steinmatrix})).  
Exploiting these we obtain
new representation
formulae for the Fisher information $J(X) = E \left[ \rho_X(X)
  \rho_X(X)^T \right]$ of an arbitrary random vector $X$
(see Theorem \ref{the:an-identity}). Such results are
  akin to those from \cite{guo2005mutual} (see Theorem \ref{the:verdu}
  for details) and, more generally, to the classical representation
  formulae for Fisher information in terms of conditional expectations
  (see \cite{LC98} or \cite{MB07} for a discussion). As an application
  we obtain (under the assumption that the Stein matrix exists)
  new information theoretic bounds for Gaussian approximation
problems. Our bounds are of the same order as those
  obtained in the pathbreaking references \cite{Jo04,BaBaNa03} (in the
  univariate setting under the assumption of a finite Poincar\'e
  constant/spectral gap; see also \cite{BaNg} for multivariate
  extensions).

Our approach is inspired by results usually exploited within the
context of the so-called \emph{Stein's method} (see
\cite{NP11,ChGoSh11}). The connexion between Stein's
  method and Fisher information was discovered in
  \cite{barbour2010compound} (in the context of compound Poisson
  approximation) and first studied explicitly by \cite{LS13b,LS12a} as
  well as \cite{sason2012entropybern}. We conclude the paper in
Section \ref{sec:conv-result-terms} with a new proof of these bounds;
our take on these matters does not rely on Stein's method and is of
independent interest to the ISIT crowd. 

The outline of the paper is as follows. All formulae and definitions
are given in Section \ref{sec:introduction}. In Section
\ref{sec:an-identity} we prove the representation
  formulae for the score in terms of the Stein matrix. In
Section~\ref{sec:conn-with-form} we prove a version of the celebrated
MMSE formula from \cite{guo2005mutual}.
In~Section~\ref{sec:conv-result-terms} we provide a general
\emph{``Stein bound''} on the standardised Fisher information of sums
of iid random vectors.

\section{Score  and Stein matrix}
\label{sec:introduction}

Fix an integer $d \ge 1$.  Let $X, Y$ be centered random $d$-vectors
(all elements in $\R^d$ are taken as $d \times 1$ column vectors)
which we throughout assume to admit a density (with respect to the
Lebesgue measure) with support $S\subset \R^d$.
\begin{definition} \label{def:score-function-stein} The \emph{score}
  of $X$ is the random vector $\rho_X(X)$ which satisfies
\begin{equation}
  \label{eq:scorefunction}
  E \left[ \rho_X(X) \varphi(X)  \right] = - E \left[ \nabla \varphi(X) \right]
\end{equation}
(with $\nabla$ the usual gradient in $\R^d$) for all test functions
$\varphi \in C_c^{\infty}(\R^d)$.
Any random  $d \times d$ matrix $\tau_X(X)$ which satisfies
\begin{equation}
  \label{eq:steinmatrix}
  E \left[ \tau_X(X) \nabla \varphi(X) \right] = E \left[ X \varphi(X) \right]
\end{equation}
for all test functions $\varphi \in C_c^{\infty}(\R^d)$ is called a
 \emph{Stein matrix} for $X$.
\end{definition}

If $X$ has covariance matrix $C$, then a direct application of the
definition of the Stein matrix  yields $E \left[ \tau_X(X)  \right] = C$;
$E \left[ \rho_X(X) \right] = 0$ and $E \left[ \rho_X(X)X^T \right] = -Id$,
where $\cdot^T$ denotes the transpose operator and $Id$ is the $d\times d$
identity matrix.
For a Gaussian random vector $Z$ with covariance matrix $C$ one uses
the well-known Stein identity (see, e.g., \cite{Jo04})
\begin{equation}\label{eq:steinidentity}
  E \left[ Z \varphi(Z)  \right] = C E \left[ \nabla \varphi(Z) \right]
\end{equation}
to prove that $ \rho_Z(Z) = -C^{-1}Z$  is the score of $Z$
and $\tau_Z(Z) = C$ is a Stein matrix of $Z$. 
Identity (\ref{eq:steinidentity}) characterizes the Gaussian
distribution in the sense that a random vector $X$ with support $\R^d$
satisfies (\ref{eq:steinidentity}) for all $\varphi \in
C_c^{\infty}(\R^d)$ if and only if $X$ is itself
Gaussian with covariance $C$. More generally,  the following result
holds (see, e.g., \cite{Jo04}).

\begin{prop}
 Let $X$ have density $f$. If $X$ has a score then it
   is uniquely defined as $\rho_X(X)$ with
 $\rho_X(x)= \nabla \log f(x)$.  
\end{prop}
\noindent

In the case $d=1$, under standard assumptions of regularity of the
density $f$, the existence of the Stein matrix $\tau$ 
  follows from standard integration by parts arguments, from
which one deduces that $\tau$ is uniquely defined as
$\tau(x) = f(x)^{-1} \int_x^\infty f(y)dy.  $ 
In higher dimensions, the existence of a Stein matrix for $X$ also
follows easily from an integration by parts argument, once one can
find a matrix valued function $x\mapsto A(x)$ whose components
$a_{ij}$ with $1 \le i,j \le d$ satisfy $\sum_{j=1}^d
\frac{\partial}{\partial x_j}\left( a_{ij}(x) f(x)\right) = -x_i$ for
all $i = 1, \ldots, d$.  As demonstrated in the huge body of
literature revolving around Malliavin calculus (see
\cite{nourdin2013entropy} as well as the monograph \cite{NP11}), a
Stein matrix always exists for random vectors that are given by a
smooth transformation of a given Gaussian field. Contrarily to the
score, however, there is no reason for which the Stein matrix, at
least according to our definition, should be unique.

 \begin{definition}
Let $X$ be a $d$-random vector with density $f$ and covariance
$B$ (invertible), and let $\phi$ be the density of a centered Gaussian
random vector $Z$ with covariance $C$.   The \emph{relative  entropy}
of $X$ is $D(X\, \| \, Z) = E \left[ \log (f(X)/\phi(X)) \right]$. The
  \emph{Fisher information (matrix)} of $X$
is $J(X) = E \left[ \rho_X(X) \rho_X(X)^T \right]$ and 
  its \emph{relative Fisher information matrix}  is
 $\mathcal{J}(X) = E
\left[ ( \rho_X(X)+B^{-1}X)( \rho_X(X)+B^{-1}X)^T 
\right]$. The
 \emph{standardised Fisher
 information distance} of $X$ is
$ J_{st}(X) = {\rm tr}\left(B \mathcal{J}(X) \right)$,
with `${\rm tr}$' the usual trace operator.
\end{definition}

Entropy and Fisher information are related to one another
via  the so-called \emph{de Bruijn's identity}, see   \cite[Lemma
2.2]{MR1851387} for the original statement, as well as
\cite[Lemma 2.3]{nourdin2013entropy} for the forthcoming version.

\begin{lemma}[Multivariate de Bruijn's identity]\label{l:db} Let  $X$
  be a random $d$-vector with covariance $C$ (invertible) and 
  let $Z$ be Gaussian  with covariance $C$ as well. Then
$D(X \, \| Z) = \int_0^1 \frac{1}{2t}  J_{st}(X_t) dt$.
 \end{lemma}
 \begin{remark}
There is some confusion  
  surrounding the denomination ``de Bruijn's identity''
  as several different (and not perfectly equivalent) formulations of this
  identity are available in the literature. See e.g. \cite[Section
  II.D]{guo2005mutual} for an alternative formulation.
 \end{remark}

\section{Representation formulae}
\label{sec:an-identity}

The following lemma is a generalization of 
\cite[Lemma 2.9]{nourdin2013entropy} to the case of summands with
arbitrary distribution. The device contained in the proof (namely a
probabilistic integration by parts formula) will be used
throughout the subsequent arguments.
\begin{lemma}\label{lem:npsformula}
  Let $X$ and $Y$ be stochastically independent centered random vectors in
  $\R^d$. Suppose that $X$ (resp.,  $Y$) has score  $\rho_X(X)$
  (resp.,   $\rho_Y(Y)$) and Stein matrix $\tau_X(X)$ (resp.,
  $\tau_Y(Y)$). For  $0 < t < 1$, let $W_t =
  \sqrt{t}X + \sqrt{1-t}Y$ and  $\Gamma_t$ be the covariance matrix of
  $W_t$. Then 
  \begin{eqnarray}
&&   \rho_{W_t}(W_t)+\Gamma_t^{-1}W_t\label{eq:pechepouletd}\\
& =&  E \left[ \frac{t}{\sqrt{1-t}} (Id - \Gamma_t^{-1}\tau_X(X))
      \rho_Y(Y)\right.\notag\\
      &&\left.+ \frac{1-t}{\sqrt{t}} (Id - \Gamma_t^{-1}\tau_Y(Y))
      \rho_X(X)\,\big| \, W_t \right] \notag
  \end{eqnarray}
is a version of the score of $W_t$. 
\end{lemma}
\noindent
{\it Proof}.
Let $\varphi \in C_c^\infty(\R^d)$ be a test function. Applying first
(\ref{eq:scorefunction}) (with respect to $Y$) then
(\ref{eq:steinmatrix}) (with respect to $X$) we get
\begin{align*}
 &  \frac{1}{\sqrt{1-t}}  E \left[  (Id -
      \Gamma_t^{-1}\tau_X(X)) \rho_Y(Y) \varphi(W_t) \right] \\
    & =  -E \left[  (Id -
      \Gamma_t^{-1}\tau_X(X)) \nabla \varphi(W_t) \right] \\
    & = - \left( E \left[ \nabla \varphi(W_t) \right]  -\Gamma_t^{-1} \frac{1}{\sqrt{t}}
    E \left[
      X   \varphi(W_t) \right] \right).
\end{align*}
Likewise
\begin{align*}
&    \frac{1}{\sqrt{t}}\lefteqn{  E \left[ E \left[ (Id -
     \Gamma_t^{-1} \tau_Y(Y)) \rho_X(X)\, | \, W_t \right]
   \varphi(W_t) \right]}\\
&   = -\left( E \left[ \nabla \varphi(W_t) \right]  -\Gamma_t^{-1} \frac{1}{\sqrt{1-t}}
    E \left[
      Y   \varphi(W_t) \right] \right).
\end{align*}
Hence
\begin{align*}
&   E \left[ E \left[ \frac{t}{\sqrt{1-t}}(Id - \Gamma_t^{-1}
      \tau_X(X)) \rho_Y(Y)\right.\right.\\
     & \left.\left. + \frac{1-t}{\sqrt{t}} (Id - \Gamma_t^{-1}
      \tau_Y(Y)) \rho_X(X)\,\big| \, W_t \right] \varphi(W_t)
  \right] \\
  & = -E \left[ \nabla \varphi(W_t) \right] + E \left[ \Gamma_t^{-1}W_t
    \varphi(W_t) \right]\\
    &=E\left[(\rho_{W_t}(W_t)+\Gamma_t^{-1}W_t)\varphi(W_t)\right],
\end{align*}
and the conclusion (\ref{eq:pechepouletd}) follows.
\qed

It is immediate to   extend (\ref{eq:pechepouletd}) to an
arbitrary number of summands. 
\begin{lemma}\label{eq:pechepouletd2} Let $X_i$, $i=1, \ldots, n$ be independent random 
  vectors with Stein matrices $\tau_i=\tau_{X_i}$ and score functions
  $\rho_i=\rho_{X_i}$, $i=1, \ldots, n$.  For all $t = (t_1, \ldots, t_n) \in [0,1]^d$
  such that $\sum_{i=1}^n t_i=1$ we define $W_t= \sum_{i=1}^n
  \sqrt{t_i}X_i$ and denote $\Gamma_t$ the corresponding covariance
  matrix. Then $\rho_t(W_t) + \Gamma_t^{-1}W_t = \sum_{i=1}^n \!
  \frac{t_i}{\sqrt{t_{i+1}}} E \left[ \left( Id -
      \Gamma_t^{-1}\tau_i(X_i) \right) \rho_{i+1}(X_{i+1}) | W_{t}
  \right]$ where we identify $X_{n+1} = X_1$ and $t_{n+1} = t_1$, and
  where we set $\rho_t=\rho_{W_t}$.
\end{lemma}

In \cite{nourdin2013entropy} we use a  version of
(\ref{eq:pechepouletd}) specialised to the case where
$X$ has covariance $C$ and  $Y=Z$ is a Gaussian random vector also with
covariance $C$. Then $\Gamma_t = C$ and, setting 
$X_t = \sqrt{t}X + \sqrt{1-t}Z$, we get, for all $0 < t < 1$,
\begin{eqnarray}
  \label{eq:10}
  &&\rho_t(X_t) + C^{-1} X_t \\
  &=& \!\!\!\!-\frac{t}{\sqrt{1-t}} E \left[ \left(
      Id- C^{-1} \tau_{X}(X) \right) C^{-1}Z \, | \, X_t\right].\notag
\end{eqnarray}
Taking
squares and simplifying accordingly we obtain the following representations for the Fisher
information and the standardised Fisher information
of an arbitrary random vector with density.
\begin{theorem}\label{the:an-identity}
Let $X$ be centered with covariance $Id$  independent of $Z$ standard
Gaussian  and $X_{t} = \sqrt{t} X + \sqrt{1-t}Z$. For all $0<t<1$, $J(X_t)$ equals
\begin{align}
  \label{eq:16}
& \frac{t^2}{1-t} E \!\!\left[  \!E \!\left[ \left(
        Id- C^{-1} \tau_{X}(X) \right) C^{-1}Z  |  X_t\right]
  \right.  \\
&  \quad\times \left.\!E \!\left[ \left(
        Id- C^{-1} \tau_{X}(X) \right) C^{-1}Z  |  X_t\right]^{T}\right]\!+ C^{-1},\nonumber
\end{align}
and 
$J_{st}(X_t)$ equals
\begin{align}
  \label{eq:11}
& = \frac{t^2}{1-t}  {\rm tr}\!\left( \!C E \!\left[ \!E \!\left[
      \left(   Id- C^{-1} \tau_{X}(X) \right) C^{-1}Z  |  X_t\right] \right.\right.
    \nonumber \\
&\quad \times \left.\left.E \!\left[
      \left(   Id- C^{-1} \tau_{X}(X) \right) C^{-1}Z  |  X_t\right] ^T  \right] \!\right). 
\end{align}
\end{theorem}
 
Arguably, the main application of formula \eqref{eq:11} provided in the present paper appears in Section \ref{sec:conv-result-terms}, where we will deduce explicit bounds in the multidimensional entropic CLT. However, representation results such as \eqref{eq:11} cover a much wider ground of applications, as they may lead in principle to new identities and new estimates in any domain where information functionals do appear. The reader is referred e.g. to \cite{LNP14} for a panoply of novel applications of formulae analogous to \eqref{eq:11} to log-Sobolev and transport inequalities.

\section{Connection with a formula of Guo, Shamai and Verd\'u}
\label{sec:conn-with-form}

It was brought to our attention (by Oliver Johnson, personal
communications) that representation (\ref{eq:10}) resembled, at least
in principle, an identity for Fisher information discovered in
\cite{guo2005mutual}. The purpose of this section is to make the  connection between the two
approaches explicit. 

\begin{lemma}[\textcolor{black}{\cite[equation (56)]{guo2005mutual}}] 
  Let $X$ be a centered random vector with covariance $C$
  independent of  $Z$  Gaussian with the same covariance
  as $X$. Then, for all $0<t<1$,  the random vector  $X_t = \sqrt{t}X
  + \sqrt{1-t}Z$ has a score
$\rho_t(X_t) = -\frac{1}{1-t}C^{-1} \left( X_t - \sqrt{t}E \left[ X \, |
        \, X_t \right] \right)$
and its Fisher information $J(X_t)$ equals
$\frac{C^{-1}}{1-t} - \frac{t\,C^{-1}}{(1-t)^2}  E \left[ \left( X\!-\!E
      \left[ X  |  X_t \right] \right) \!\left( X\!\!-\!E
      \left[ X  |  X_t \right] \right)\!^T\right]\!C^{-1}.$
\end{lemma}
\noindent
{\it Proof.}
Clearly,  $X_t$ has a differentiable density with support $\R^d$.
 Let $\varphi \in C_c^{\infty}(\R^d)$ be a test
  function. Then 
  \begin{align*}
  &  E \left[ C^{-1} \left( X_t - \sqrt{t}\,E \left[ X \, |
        \, X_t \right] \right) \varphi(X_t)  \right]\\
        &  =
 E \left[ C^{-1} \left( X_t - \sqrt{t} X  \right) \varphi(X_t)  \right] \\
& =\sqrt{1-t} E \left[ C^{-1}  Z\varphi(X_t)
\right]  = (1-t) E \left[ \nabla \varphi(X_t) \right].
  \end{align*}
Both   claims then follow after
  straightforward computations.
\qed

Next, as in \cite{guo2005mutual}, we define 
${\rm MMSE}(X, t)  = E \left[ \left( X-E
      \left[ X \, | \, X_t \right] \right) \left( X-E
      \left[ X \, | \, X_t \right] \right)^T\right]$.
 Direct application of the above yields  the following. 
\begin{prop}\label{the:verdu} If $A$ is a matrix we write $A^2$ for $A A^T$. Then 
  $Id - \frac{1}{1-t} {\rm MMSE}(X,t)C^{-1} = t \,C E \!\left[ \!E
    \left[ \left( Id- C^{-1} \tau_{X}(X) \right) C^{-1}Z |
      X_t\right]^2 \right]$ so that $J_{st}(X_t)= \frac{t}{1-t} {\rm
    tr}\left( Id - \frac{1}{(1-t)} {\rm MMSE}(X,t)C^{-1} \right)$.
\end{prop}
Plugging this last identity into  Lemma 
 \ref{l:db} shows that relative entropy is an integral of minimal
 squared error; this claim is equivalent (up to scaling) to    \cite[equation
 (57)]{guo2005mutual}.

\section{Information bounds for sums of random vectors}
\label{sec:conv-result-terms}

In the sequel we suppose for simplicity that all random vectors 
are isotropic (i.e. have  identity covariance matrix). 
 \begin{theorem}\label{thm:20}
Let $X_1, \ldots, X_n$ be independent random vectors in $\R^d$ and suppose
that the $X_i$ have Stein matrix $\tau_i(X_i)$ and score function
$\rho_i(X_i)$. 
   Let $W_n = \frac{1}{\sqrt{n}} \sum_{i=1}^nX_i$. Define $W_n^{(t)} =
 \sqrt{t}W_n + \sqrt{1-t}Z$, where $Z$ is an independent standard Gaussian random vector. Then
$
{J}_{st}(W_n^{(t)})   \le \frac{t^2}{n^2(1-t)} \sum_{i=1}^n
 \!  {\rm tr}\! \left( \! E \! \left[ \! \left( Id-\tau_i(X_i) \right)\! \left( Id-\tau_i(X_i) \right)^T \! \right] \! \right)$
for all $0 \le t \le 1$.
\end{theorem}
\noindent
{\it Proof}.
 First, by Jensen's
    inequality, we see that
      \begin{align*}
&        {\rm tr} \big(E [ E [ (Id- 
 \tau_{W_n}(W_n))Z \, | \, W_n^{(t)}] \\
& \quad \times E [ (Id- 
 \tau_{W_n}(W_n))Z \, | \, W_n^{(t)}]^T]\big) \\
& \le {\rm tr}\big( E [  (Id- 
 \tau_{W_n}(W_n))  (Id- 
 \tau_{W_n}(W_n))^T]\big).
      \end{align*}
Next, it is easy to prove (see  \cite{nourdin2013entropy} for a proof when
$d=1$) that  $\tau_{W_n}(W_n) = \frac{1}{n} \sum_{i=1}^n E \left[ \tau_i(X_i) \, |
      \, W_n \right]$ is a Stein matrix for $W_n$.  
Hence, by (\ref{eq:11}), $J_{st}(W_n^{(t)})$ is less than or equal to
\begin{align*}
&   \frac{t^2}{1-t} {\rm tr}\left( E \left[ (Id -
    \tau_{W_n}(W_n))( Id - \tau_{W_n}(W_n))^T \right] \right) \\
& \le \frac{1}{n^2}  \frac{t^2}{1-t} {\rm tr} \left(  E \left[ \left(
      \sum_{i=1}^n ( Id - \tau_i(X_i)) \right)\right.\right.\\
& \quad \quad \quad \times \left.\left.\left( \sum_{i=1}^n ( Id - \tau_i(X_i))^T \right)\right]
 \right).
\end{align*}
Independence of the $X_i$ as well as the fact that $ E \left[ Id -
\tau_i(X_i)  \right]=0$ allow to  conclude. 
\qed

In particular, if the $X_i$ are i.i.d. copies of $X$ then
${J}_{st}(W_n^{(1/2)})   \le \frac{1}{2n}
  {\rm tr}\left(  E \left[ \left( Id-\tau_X(X) \right)\left(
        Id-\tau_X(X) \right)^T \right] \right).$ By Cramer's theorem (see, e.g., \cite{kagan1989multivariate}),
convergence of $W_n$ to the Gaussian is equivalent to convergence of
$W_n^{(1/2)}$, and Theorem \ref{thm:20} provides rates of convergence
(of order $1/n$) of the Fisher information under the assumption that
$X$ has a well-defined Stein matrix $\tau_X(X)$. 
A straightforward extension of  \cite[Lemma 1.21]{Jo04}
   to the multivariate setting 
   shows that standardised information decreases along
 convolutions.
\begin{lemma}\label{eq:35} If $X$ and $Y$ are independent
  isotropic (i.e.  identity covariance matrix) real-valued $d$-random
  vectors then $J_{st}(\sqrt{t} X+ \sqrt{1-t}Y) \le t
  J_{st}(X) + (1-t) J_{st}(Y)$.
\end{lemma}
\noindent 
{\it Proof}. Let $W_t=\sqrt{t} X+ \sqrt{1-t}Y$. From
definition \eqref{eq:scorefunction} it is easy to see that $\rho_t(w)
= E \left[ \sqrt{t} \rho_X(X)+ \sqrt{1-t}\rho_Y(Y) \, | \, W_t=w
\right]$. By definition $J(W_t) =
\sum\limits_{j=1}^dE \left[(\rho_t(W_t))_j^2\right]$ which, by Jensen's
inequality, is smaller or equal to $
 t E \left[\sum\limits_{j=1}^d (\rho_X(X))_j^2 \right] + (1-t)  E
 \left[\sum\limits_{j=1}^d (\rho_Y(Y))_j^2 \right]  = tJ(X) + (1-t)
 J(Y)$, and the claim is proved for Fisher information. The extension
 to $J_{st}$ is immediate. \qed 

In particular, from Lemma \ref{eq:35}, if 
$Z$ is standard Gaussian independent of $X$, then 
$J_{st}(X_t) \le t J_{st}(X) + (1-t) J_{st}(Z) = t J_{st}(X)$
for all $0 \le t \le
1$ with $X_t = \sqrt{t} X + \sqrt{1-t}Z$, so that
\begin{equation}
  \label{eq:33}
D(X\, \| \, Z) \le \frac{1}{2} J_{st}(X)
\end{equation}
by Lemma \ref{l:db}. 
Hence  bounds  on the standardised Fisher information
translate directly into bounds on the relative entropy hereby
providing, via Pinsker's inequality
\begin{equation}\label{pinsker}
2  d_{TV}(X, Z) \le \sqrt{ 2D(X \, \| \, Z)},
\end{equation} 
bounds on the \emph{total variation distance}
between the law of $X$ and the law of $Z$.  From
  \eqref{eq:33} we  thus obtain rates of convergence in total
variation which have the correct order (see
e.g. \cite{MR2128239,BaBaNa03} for similar rates of convergence under
the assumption of finite Poincar\'e constant).

\begin{remark} {\rm 
It is
      still largely an open question how the assumption of existence
      of a Stein matrix relates with more
      standard assumptions such as  finiteness of the  Poincar\'e
      constant. See  e.g. \cite{airault2010stokes, nourdin2013entropy, LNP14} for 
      discussions.  
  }
\end{remark}

\section{Stein representations for Fisher information}
\label{sec:stein-repr-fish}

Our next lemma
provides a new handle on conditional expectations
  which is also of independent interest.
\begin{lemma}[Poly's lemma]\label{lem:poly}
Let $X$ and $Y$ be square-integrable random variables with mean $E[X]=0$. Then
$E \left[ \left(  E \left[  X \, |\, Y \right]     \right)^2  \right] =
 \sup_{\varphi\in\mathcal{H}(Y)}  \left( E \left[ X \varphi(Y) \right] \right)^2$,
where the supremum is taken over the collection $\mathcal{H}(Y)$ of  functions $\varphi$ such
that $E [\varphi(Y)] = 0$ and $ E \left[   \varphi(Y)^2 \right]\le1$.
\end{lemma}
\noindent
{\it Proof}.
First, by
  Cauchy-Schwarz,
  \begin{align*}
&    \sup_{\varphi\in\mathcal{H}(Y)} \left(E  \left[ X \varphi(Y) \right]\right)^2 =     \sup_{\varphi\in\mathcal{H}(Y)} \left(E  \left[ E [X|Y]
      \varphi(Y) \right]\right)^2 \\
& \le \!\!\sup_{\varphi\in\mathcal{H}(Y)}\!\! E \left[ E [X | Y]^2 \right] \!E
    \left[ \varphi(Y)^2 \right] \!\le \!E \!\left[ E [X | Y]^2 \right]\!.
  \end{align*}
To prove the reverse inequality define  $\varphi(y) = E \left[ X | Y=y \right]/\sqrt{E \left[ E [X
    | Y]^2 \right]}$. Clearly $E [\varphi(Y) ]= 0$ and $E [\varphi(Y)^2] \le
1$  so that $\varphi \in \mathcal{H}(Y)$ and $\sup_{\varphi\in\mathcal{H}(Y)} \left( E \left[ X \varphi(Y) \right]\right)^2$ is bigger than or equal to
\begin{align*}
    & \left(E \left[ X \frac{E \left[ X | Y \right]}{\sqrt{E \left[ E [X
    | Y]^2 \right]}} \right]\right)^2 \\
& = \frac{(E \left[ X  E \left[ X | Y \right] \right])^2}{ E \left[ E [X
    | Y]^2 \right]} =  E \left[ E \left[ X | Y \right]^2 \right].
\end{align*}
Equality ensues.
\qed

We immediately deduce an original proof (not relying on Stein's method!) of a
recently discovered fact (see e.g.
\cite{barbour2010compound,LS13b,LS12a}) that the Fisher
information distance is dominated by expressions which appear
naturally within the context of Stein's method.

\begin{theorem}[Stein representation for relative Fisher information, $d=1$]\label{sec:stein-repr-fish-1}
  Let $W_n = \frac{1}{\sqrt{n}} \sum_{i=1}^nX_i$ where the $X_i$ are
  independent random variables with Stein factor $\tau_i(X_i)$ and
  score function $\rho_i(X_i)$. Then $ \mathcal{J}(W_n) =
  \sup_{\varphi\in\mathcal{H}(W_n)} \left( E \left[ \varphi'(W_n) -
      W_n \varphi(W_n) \right] \right)^2.$
\end{theorem}
\noindent
{\it Proof}.
  We combine Lemma \ref{eq:pechepouletd2} (in the
  special case $t_1 = t_2 = \ldots = t_n = 1/n$) and Lemma \ref{lem:poly}
  to deduce 
  that $n\mathcal{J}(W_n)$ is given by
  \begin{align*}
 &\quad\,\,  n\,E \left[ (\rho_n(W_n)+W_n)^2 \right]  \\
& =  E\big[ \big(  E[\sum_{i=1}^n
          (1-\tau_i(X_i))\rho_{i+1}(X_{i+1})   |  W_n ] \big)^2\big] \\
  & = \!\!\!\!\!\sup_{\varphi\in\mathcal{H}(W_n)} \!\!\!\!\big( \!\sum_{i=1}^n
  \!\!E \!\left[
          (1-\tau_i(X_i))\rho_{i+1}(X_{i+1}) \varphi(W_n)  \right]\!\big)^2 \\
      & =  \sup_{\varphi\in\mathcal{H}(W_n)} E \left[ \sqrt{n}
         \varphi'(W_n) - \sum_{i=1}^n X_i \varphi(W_n) \right]^2,
  \end{align*}
and the conclusion follows.
\qed

\begin{cor}[Stein representation for relative Fisher information]
  Let $W_n = \frac{1}{\sqrt{n}} \sum_{i=1}^nX_i = ((W_n)^1, \ldots,
  (W_n)^d)^T$ where the $X_i$ are independent $d$-random vectors with
  Stein matrix $\tau_i(X_i)$ and score function $\rho_i(X_i)$. Then $
  J_{st}(W_n) = \sum\limits_{j=1}^d\sup_{\varphi\in\mathcal{H}(W_n)}
  \left( E \left[
        \partial_j \varphi(W_n) - (W_n)^j \varphi(W_n) \right] \right)^2.$
\end{cor}

\section*{Acknowledgment}

We are grateful to Oliver Johnson for providing the connexion with
\cite{guo2005mutual} and to Guillaume Poly for sharing his Lemma
\ref{lem:poly} with us.  We thank the referees for
  their careful comments which helped improve the rendition of this
  work. Giovanni Peccati was partially supported by the Grant F1R-MTH-PUL-12PAMP  (PAMPAS) from Luxembourg University.  Yvik Swan gratefully acknowledges support from the IAP
  Research Network P7/06 of the Belgian State (Belgian Science
  Policy).

\end{document}